\documentclass[12pt]{article}
\usepackage{amsmath, graphicx}
\usepackage{amstext,amsfonts,amsbsy,eucal,amssymb, amsthm, color}

\newtheorem{theorem}{Theorem}[section]

\theoremstyle{definition}

\newtheorem{example}[theorem]{Example}
\newtheorem{remark}[theorem]{Remark}

\def\Z{\operatorname{\mathbb{Z}}}

\def\R{\operatorname{\mathbb{R}}}
\def\C{\operatorname{\mathbb{C}}}
\def\P{\operatorname{\mathbb{P}}}

\def\pic{\operatorname{{\rm Pic}}}
\def\H{\operatorname{{\mathcal H}}}
\def\E{\operatorname{{\mathcal E}}}
\def\F{\operatorname{{\mathcal F}}}

\begin{document}
\begin{center}
{\large A note on minimization of rational surfaces obtained from birational dynamical systems} 
\end{center}
\begin{center}
Adrian Stefan Carstea, Tomoyuki Takenawa
\end{center}
\begin{center}
{\it National Institute of Physics and Nuclear Engineering, Dept. of Theoretical Physics, Atomistilor 407, 077125, Magurele, Bucharest, Romania}\\ 
{\it Faculty of Marine Technology, Tokyo University of Marine Science and Technology, 2-1-6 Etchu-jima, Koto-Ku, Tokyo, 135-8533, Japan}\\ 
\end{center}

\begin{abstract}
In many cases rational surfaces obtained by desingularization of birational dynamical systems are not relatively minimal.
We propose a method to obtain coordinates of relatively minimal rational surfaces by using blowing down structure. We apply this method 
to the study of various integrable or linearizable mappings, including discrete versions of reduced Nahm equations. 
\end{abstract}

\section{Introduction}

Studies of (possibly non-autonomous) discrete integrable dynamical 
systems started with discovery of singularity confinement method \cite{singconf}, which imposes conditions on the singularities of solutions, and many discrete analogues of Painlev\'e equations were found through this method. In \cite{Sakai01} Sakai, noting Okamoto's earlier work on continuous case \cite{okamoto}, showed that all the continuous or discrete Painlev\'e equations are classified in terms of so called generalized Halphen surfaces, where the condition for the dimension of the linear system of the anti-canonical divisor is relaxed to be not necessarily one. In these studies surfaces depend on parameters including the independent variable, and each surface is called the space of initial conditions.  In discrete case, the dynamical system is lifted to a sequence of isomorphisms between surfaces.
  
On the other hand, in \cite{DF01} Diller and Favre showed that for any birational automophism $\varphi$ on a projective smooth rational surface $S$, we can construct a rational surface $\tilde{S}$ by successive blow-ups from $S$ such that (i) $\sigma\circ \tilde{\varphi}= \varphi \circ\sigma$, where $\sigma$ denotes the successive blow-downs $\sigma: \tilde{S}\to S$, (ii) $\tilde{\varphi}=\varphi$ on $S$ in generic, and (iii) $\tilde{\varphi}:\tilde{S}\to \tilde{S}$ is analytically stable.
In general, $\tilde{\varphi}$ is said to be lifted from $\varphi$ 
if the condition (i) and (ii) are satisfied and
a birational automorphism  $\varphi$ on $S$ is said to be analytically stable if the condition $((\varphi)^*)^n=((\varphi)^n)^*$ holds on the Picard group on $S$. The notion of analytical stability is closely related to the singularity confinement. Indeed, this notion is equivalent to the condition that
there is no curve $C$ on $S$ and a positive integer $k$ such that 
$\varphi(C)$ is a point on $S$ and $\varphi^k(C)$ is an indeterminate point of $\varphi$, i.e. analytical stability demands that singularities are not recovered by the dynamical system. 
In order to compare these two notions, let us consider a mapping $\varphi$ having both confined and unconfined sequences. For a confined sequence there is a curve such that $\varphi^i(C)$ ($1\leq i \leq k$) is a point and 
$\varphi^{k+1}(C)$ is a curve again for some positive integer $k$. If we blow up the phase space at $\varphi^i(C)$
for $1\leq i\leq k$, then the singularity would be relaxed and resolved by
successive applications of this procedure. By applying this to all
the confined sequences, we would obtain a surface
where the lifted birational automorphism $\tilde{\varphi}$ is analytically stable.

From Diller and Favre's work, such birational automorphisms $f$
are classified as follows:
Let $f$ be a bimeromorphic automophism of a K\"ahler surface with the
maximum eigenvalue of $f^*$ is one. Up to
bimeromorphic conjugacy, exactly one of the following holds.
\begin{itemize}
\item The sequence $||(f^n)^*||$ is bounded, and $f^n$ is an automorphism isotopic to the identity for some $n$, where $||\cdot ||$ denotes the Euclidean norm
w.r.t. some basis of the Picard group. 
\item The sequence $||(f^n)^*||$ grows linearly, and $f$ preserves a rational fibration. In this case, $f$ can not be conjugated to an automorphism.
(We say $f$ is linearizable or linearizable in cascade in this case \cite{ROG00,TEGORS}). 
\item The sequence $||(f^n)^*||$ grows quadratically, and $f$ is an automorphism preserving an elliptic fibration.
\end{itemize}

Though in this paper we consider mainly autonomous case, 
the procedure constructing analytically stable mapping can be applied also for non-autonomous case such as linearizable mappings or discrete Painlev\'e equations. In this case, we start from a sequence of birational mappings 
$\varphi_i: S\to S$ and blow up successively at confined singular points whose positions depend on $i$. Then, $\varphi_i$ would be lifted to a sequence of birational mappings 
$\tilde{\varphi}_i : \tilde{S}_i\to \tilde{S}_{i+1}$ such that
(i) $\sigma_{i+1} \circ \tilde{\varphi}_i= \varphi_i \circ\sigma_i$, (ii) 
$\tilde{\varphi}_i=\varphi_i$ on $S$ in generic, and (iii) $\{\tilde{\varphi}:\tilde{S}_i\to \tilde{S}_{i+1}\}_{i\in \Z}$ is analytically stable, i.e. $\varphi_i^* \circ \cdots \circ \varphi_{i+n}^*=(\varphi_{i+n} \circ \cdots \circ \varphi_i)^*$ holds on the Picard group on 
$\tilde{S}_{i+n+1}$ for any $i$ and non-negative integer $n$
(see \cite{Takenawa1, Takenawa2,TEGORS} about computation and the relation
to the degree growth).

In the study of integrable systems,
we often want to find conserved quantities or linearize a given integrable mappings, but the above construction of analytically stable mapping
does not guarantees that\\
(i) $\tilde{\varphi}$ is an automorphism.\\
(ii) $\tilde{S}$ is relatively minimal, i.e. there does not exists
a blow-down of $\tilde{S}$, $\pi:\tilde{S}\to \tilde{S}'$ such that
$\tilde{\varphi}'$ is still analytically stable on $\tilde{S}'$.\\
In other words, the following possibilities remain.\\
(ia) A singularity sequence consists of infinite sequences of points to both 
sides and a finite sequence of curves: 
$$\cdots \to {\rm point}\to {\rm point}\to {\rm curves}
\to \cdots \to{\rm curves} \to {\rm point}\to {\rm point}\to \cdots,$$
where the image of a curve $C$, parametrized as $(f(t),g(t))$ on some coordinates, under $\varphi^n$ is defined as 
the Zariski closure of 
$\lim_{\varepsilon\to 0} \varphi^n (f(t)+c_1\epsilon,g(t)+c_1\epsilon)$
 with generic $t$, $c_1$ and $c_2$.\\
(ib) A singularity sequence consists of a semi-infinite sequence of points 
and a semi-infinite sequence of curves to each side:
$$\cdots \to {\rm point}\to {\rm point}\to  \cdots \cdots \to {\rm curves}\to {\rm curves}\to \cdots.$$
(ii)' A finite set of exceptional curves are permuted.\\
 
Proposition 1.7 and Lemma 4.2 of \cite{DF01} (cf. \cite{Gizatullin80}) says that curves in (ia) can be blown down,
and that (ib) occurs only if $f$ is not conjugate to an automorphism, i.e.
if $f$ is linearizable or has a positive entropy. 
And if a curve in Case (ii)' or Case (ib) is exceptional of the first kind, we can blow down them. Hence theoretically, we can obtain relatively minimal analytically stable surfaces and compute the action on the Picard group. 
However, for investigating
properties of the mapping $f$ such as conserved quantities,
we need to know coordinate change explicitly.

Our aim in this paper
is to develop a method to control blowing down structures on the level of coordinates.
We apply our method to various examples, including the newly studied discretization of reduced Nahm equations \cite{PPS11}. 
In general, finding elliptic fibration for an elliptic surface 
is not easy if the surface is not minimal, and
we use information of singularity patterns of the dynamical systems
for finding unnecessary $(-1)$ curves. Accordingly, this method of minimization shows how an order-two mapping with complicated singularity structure can be brought to a simpler form which enables computations of conserved quantities.

In the next section, we recall some basic notions and
blowing down structures. In Section 3, we investigate
discrete versions of reduced Nahm equations, which 
 preserve a rational elliptic fibration.
We will show that the associated surfaces are not minimal and by minimization one can transform the mappings to simpler ones.
In Section 4, we investigate linearizable dynamical systems,
including non-autonomous case.

\section{Blowing down structure}

\subsection{Preliminaries}
Notations: cf. \cite{hartshorne, IS96}
\begin{align*}
S:&\ \mbox{a smooth rational surface}\\
{\mathcal D}:&\ \mbox{the linear equivalent class of a divisor $D$}\\
D\cdot D':&\ \mbox{the intersection number of divisors $D$ and $D'$} \\
{\mathcal O}(D):&\ \mbox{the invertible sheaf corresponding to $D$}\\
\pic(S)=&\ \mbox{the group of isomorphism classes of invertible sheaves on $S$}\\
\simeq &\ \mbox{the group of linear equivalent classes of divisors on $S$}\\
{\mathcal E}: &\ \mbox{the total transform of divisor class of a line on $\P^2$} \\ 
{\mathcal H}_x, {\mathcal H}_y:&\ \mbox{the total transform of divisor class of a line $x={\rm constant}$}\\ &\ \mbox{(or $y={\rm constant}$) on $\P^1\times \P^1$}\\
{\mathcal E}_i:&\ \mbox{the total transform of the exceptional divisor class of the $i$-th blow-up} \\ 
|D|\simeq &(H^0(S,{\mathcal O}(D))-\{0\})/\C^{\times} :\ \mbox{the linear system of $D$}\\
K_S: &\ \mbox{the canonical divisor of a surface $S$}\\
g(C):&\  \mbox{the genus of an irreducible curve $C$, given by the genus formula}\\&\mbox{ $g(C)=1+\frac{1}{2}(C^2+C\cdot K_S)$ if $C$ is smooth}. 
\end{align*}

{\it Blowing up:} Let $X$ be a smooth projective surface and let $p$ be a point on X. 
There exist a smooth projective surface $X'$ and a morphism $\pi: X'\to X$ such that 
$\pi^{-1}(p)\cong\P^{1}$ and $\pi$ represents a biholomorphic mapping from 
$X'-\pi^{-1}(p)\to X-(p)$. The morphism is called blowing down and the correspondence 
$\pi^{-1}$ is called blowing up of $X$ at $p$ as a rational mapping. For example if $X$ 
is the space $\C^2$ and $p$ is a point of coordinate $(x_0,y_0)$ then we denote blowing up of $X$ in
$p$ 
$$X'=\{(x-x_0,y-y_0;\zeta_0:\zeta_1)\in \C^2\times \P^1|(x-x_0)\zeta_0=(y-y0)\zeta_1\}$$
by
$$\pi:(x,y)\longleftarrow(x-x_0,(y-y_0)/(x-x_0))\cup ((x-x_0)/(y-y_0),y-y_0).$$

{\it Total transform and proper transform:} Let $\pi :Y\to X$ be the blow-down to a point $p$ on $X$ and $D$ be a divisor on $X$. 
The divisor $\pi^*(D)$ on $Y$ ($\pi^*$ denotes the pull-back from $\pic(X)$ to $\pic(Y)$) is called the total transform of $D$ and for any analytic subvariety $V$ on $X$ the closure of $\pi^{-1}(V-p)$ in $Y$ is called the proper transform of $V$.

Let $X$ be a surface obtained by $N$ times blowing up of $\P^2$. Then the Picard group Pic($X$) is isomorphic to the $\Z$-module (the Neron-Severi lattice):
$$\pic(X)=\Z {\mathcal E} \oplus \bigoplus_{i=1}^N\Z {\mathcal E}_i$$
and the intersection of two divisors on $X$ are given by the following basic formulas (valid for any $i,j=1,\cdots,N$):
$${\mathcal E}^2=1,\ {\mathcal E}_i^2=-1,\ {\mathcal E}\cdot {\mathcal E}_i=
{\mathcal E}_i \cdot {\mathcal E}_j=0\ (i\neq j).$$
The anti-canonical divisor class is
$$-K_S=3{\mathcal E} -\sum_{i=1}^N{\mathcal E}_i.$$

In the case where $X$ is a surface obtained by $N$ times blowing up of $\P^1\times \P^1$, the Picard group Pic($X$) is  
$$\pic(X)=\Z {\mathcal H}_x\oplus\Z {\mathcal H}_y\oplus \bigoplus_{i=1}^N\Z {\mathcal E}_i$$
and the intersection of divisors and the anti-canonical divisor are given by
$${\mathcal H}_x\cdot {\mathcal H}_y=1,\ {\mathcal E}_i^2=-1,\ {\mathcal E}_i\cdot {\mathcal E}_j={\mathcal E}_i\cdot {\mathcal H}_x={\mathcal E}_i\cdot {\mathcal H}_y={\mathcal H}_x^2={\mathcal H}_y^2=0\ (i\neq j),$$
$$-K_S=2{\mathcal H}_x+2{\mathcal H}_y -\sum_{i=1}^N{\mathcal E}_i.$$

{\it A rational elliptic surface}: A rational surface $X$ is called a rational elliptic surface if there exists a fibration given by the morphism: $\pi:X\rightarrow \P^1$ such that:
\begin{itemize}
\item for all but finitely many points $k\in {\P}^1$ the fibre $\pi^{-1}(k)$ is an elliptic curve.
\item $\pi$ is not birational to the projection: $E\times \P^1 \to\P^1$ for any curve $E$.
\end{itemize}
If no fibers contain exceptional curves of the first kind, 
the surface is called
minimal rational elliptic surface. 
Minimality is equivalent to the condition that the rank of the Picard group
is $10$.
In the non-autonomous case, 
the role of minimal rational elliptic surfaces are replaced by
``generalized Halphen surfaces'' \cite{Sakai01}.

Let $S=S_m$ be a surface obtained by successive $m$ times
blowing up from
$\P^2$ (or any rational surface) at indeterminate or extremal point of $\varphi$,
i.e. the Jacobian $\partial(\bar{x},\bar{y})/\partial(x,y)$ in some local coordinates is zero, such that $\tilde{\varphi}$ on $\S$ is analytically stable.  
Let $F_m$ be a curve on $S$ with self-intersection $-1$ and ${\mathcal F}_m$ 
be the
corresponding divisor class. Our strategy to write the blow-down 
$S_m$ along $F_m$ by coordinates is as follows.

Take a divisor class ${\mathcal F}$ such that
there exists a {\it blowing down structure} (this terminology is due to \cite{DO88}):
$S=S_m\to S_{m-1}\to S_{m-2} \to \cdots \to S_1\to \P^2$,
where $S_m \to S_{m-1}$ is a blow-down along $F_m$ and 
each $S_i\to S_{i-1}$ is a blow-down along an irreducible curve, 
such that the divisor class of lines in $\P^2$ is ${\mathcal F}$. 
Let $|{\mathcal F}|=\alpha_0f_0+\alpha_1f_1+\alpha_2f_2=0$.
Then $(f_0:f_1:f_2)$ gives $\P^2$ coordinates.

In order to find such ${\mathcal F}$ we note the following facts.

It is necessary for the existence of such a blow-down structure
that there exists a set of divisor classes ${\mathcal F}_1, \dots, {\mathcal F}_m$ such that 
$$ {\mathcal F}^2=1\quad {\mathcal F}_i^2=-1,\quad 
{\mathcal F}_i\cdot {\mathcal F}_j=0, \quad   
{\mathcal F} \cdot {\mathcal F}_i=0$$
for ($1\leq i,j \leq m$),
and further
that (i) the genus of divisor $F$ is zero;
(ii) the linear system of ${\mathcal F}$ does not have a fixed part
in the sense of Zariski decomposition and its dimension is two.

If the linear system of ${\mathcal F}$ does not have fixed part, 
then by Bertini theorem, 
its generic divisor is smooth and irreducible (this follows from the fact that two divisors defines a pencil by blowing up at the unique intersection and 
P. 137 of \cite{IS96}), and its genus is given by the formula
$$g=1+\frac{1}{2}(F^2+F\cdot K_S).$$ 
From this fact and Condition (ii), $1+\frac{1}{2}(F^2+F\cdot K_S)$
should be zero. 

\begin{example}
If degree of ${\mathcal F}$ is less than 6, then ${\mathcal F}$ is given by one of the following forms.
\begin{align}\label{classF}
&\begin{array}{ll}
{\mathcal E}\\
2{\mathcal E}-{\mathcal E}_{i_1}-{\mathcal E}_{i_2}-{\mathcal E}_{i_3}\\
3{\mathcal E}-2{\mathcal E}_{i_1}-{\mathcal E}_{i_2}-{\mathcal E}_{i_3}-{\mathcal E}_{i_4}-{\mathcal E}_{i_5}\\
4{\mathcal E}-2{\mathcal E}_{i_1}-2{\mathcal E}_{i_2}-2{\mathcal E}_{i_3}-{\mathcal E}_{i_4}-{\mathcal E}_{i_5}-{\mathcal E}_{i_6}\\
4{\mathcal E}-3{\mathcal E}_{i_1}-{\mathcal E}_{i_2}-{\mathcal E}_{i_3}-{\mathcal E}_{i_4}-{\mathcal E}_{i_5}-{\mathcal E}_{i_6}-{\mathcal E}_{i_7}\\
5{\mathcal E}-2{\mathcal E}_{i_1}-2{\mathcal E}_{i_2}-2{\mathcal E}_{i_3}-2{\mathcal E}_{i_4}-2{\mathcal E}_{i_5}-2{\mathcal E}_{i_6}\\
5{\mathcal E}-3{\mathcal E}_{i_1}-2{\mathcal E}_{i_2}-2{\mathcal E}_{i_3}-2{\mathcal E}_{i_4}-{\mathcal E}_{i_5}-{\mathcal E}_{i_6}-{\mathcal E}_{i_7}\\
5{\mathcal E}-4{\mathcal E}_{i_1}-{\mathcal E}_{i_2}-{\mathcal E}_{i_3}-{\mathcal E}_{i_4}-{\mathcal E}_{i_5}-{\mathcal E}_{i_6}-{\mathcal E}_{i_7}-{\mathcal E}_{i_8}-{\mathcal E}_{i_9},
\end{array}
\end{align}  
where $i_j$'s are all distinct with each other. All the above ${\mathcal F}$ admit blow-down structure if the positions of blow-up points are generic. For example,
for ${\mathcal F}=2{\mathcal E}-{\mathcal E}_{i_1}-{\mathcal E}_{i_2}-{\mathcal E}_{i_3}$, ${\mathcal F}_i$'s are given by
$${\mathcal E}-{\mathcal E}_{i}-{\mathcal E}_{j}
(\{i,j| i\neq j\} \subset  \{i_1,i_2,i_3\}),\
{\mathcal E}_j (j\neq i_1,i_2,i_3) $$ 
and for 
${\mathcal F}=3{\mathcal E}-2{\mathcal E}_{i_1}-{\mathcal E}_{i_2}-{\mathcal E}_{i_3}-{\mathcal E}_{i_4}-{\mathcal E}_{i_5}$,
${\mathcal F}_i$'s are given by 
$${\mathcal E}-{\mathcal E}_{i_1}-{\mathcal E}_{j} 
(j\in  \{i_2,\dots,i_5\}),\
2{\mathcal E}-{\mathcal E}_{i_1}-{\mathcal E}_{i_2}-{\mathcal E}_{i_3}-{\mathcal E}_{i_4}-{\mathcal E}_{i_5},\
{\mathcal E}_j (j\neq i_1,\dots,i_5).$$ 

In the above claim, the condition ``the positions of blow up points are generic'' is not easy to describe explicitly. For example, the dimension of the linear system of $2{\mathcal E}-{\mathcal E}_{i_1}-{\mathcal E}_{i_2}-{\mathcal E}_{i_3}$ is less than two if the base points of $i_j$-th blow ups ($j=1,2,3,4$) are on the same line, if the base point of $i_1$-th blow up is on the $i_4$-th exceptional curve, or if the base points of $i_1$-th and $i_2$-th blow ups are on the $i_3$-th exceptional curve (in the third case, the quadratic curve with the divisor class $2{\mathcal E}-{\mathcal E}_{i_1}-{\mathcal E}_{i_2}-{\mathcal E}_{i_3}$ is unique). Nevertheless, the above list is useful for finding blow down structure, since we can easily compute the linear system of a given divisor class for explicit surfaces.
\end{example}

If we want to blow down to $\P^1\times \P^1$ instead of $\P^2$,
our strategy becomes as follows.

Let $F_{m-1}$ be a curve on $S$ with self-intersection $-1$ and 
${\mathcal F}_{m-1}$ 
be the
corresponding divisor class.
Take a divisor class ${\mathcal H}_u$ and ${\mathcal H}_v$ such that
there exists a blow-down structure:
$S=S_{m-1}\to S_{m-2}\to \to \cdots \to S_1\to \P^1 \times \P^1$, 
where
$S_{m-1} \to S_{m-2}$ is a blow-down along $F_{m-1}$ and 
each $S_i\to S_{i-1}$ is a blow-down along an irreducible curve, 
such that the divisor class of lines $u={\rm const}$ and $v={\rm const}$ are
${\mathcal H}_u$ and ${\mathcal H}_v$. 
Let $|{\mathcal H}_u|=\alpha_0f_0+\alpha_1f_1=0$ and 
$|{\mathcal H}_v|=\beta_0g_0+\beta_1g_1=0$.
Then $(u,v)=(f_0/f_1,g_0/g_1)$ gives $\P^1\times \P^1$ coordinates.

In this case, it is necessary 
that there exits a set of divisor classes ${\mathcal F}_1, \dots, {\mathcal F}_{m-1}$ such that 
\begin{align*}
{\mathcal H}_u^2={\mathcal H}_v^2=0,  {\mathcal H}_u\cdot{\mathcal H}_v=1,
\quad {\mathcal F}_i^2=-1,\\ 
{\mathcal F}_i\cdot {\mathcal F}_j=0, \quad   
{\mathcal H}_u \cdot {\mathcal F}_i={\mathcal H}_v \cdot {\mathcal F}_i=0
\end{align*}
for ($1\leq i\neq j \leq m-1$),
and further
that (i) each genus of divisor $\H_u$ or $\H_v$ is zero;
(ii) each linear system of ${\mathcal H}_u$ or ${\mathcal H}_v$ does not have 
a fixed part and its dimension is one.
Consequently, $1+\frac{1}{2}(F^2+F\cdot K_S)$
should be zero again. 

\begin{example}
If $S$ is obtained by successive blow-ups from $\P^2$,
and the sum of degree of ${\mathcal H}_u$ or ${\mathcal H}_v$ is less than 6, then each ${\mathcal H}_u$ or ${\mathcal H}_v$ is given by 
${\mathcal F} -{\mathcal E}_k$, where ${\mathcal F}$ is in the list \eqref{classF}.

If $S$ is obtained by successive blow-ups from $\P^1\times \P^1$,
each ${\mathcal H}_u$ or ${\mathcal H}_v$ is given by
\begin{align}
&\begin{array}{ll}
{\mathcal H}_x \\
{\mathcal H}_x+{\mathcal H}_y-{\mathcal E}_{i_1}-{\mathcal E}_{i_2}\\
2{\mathcal H}_x+{\mathcal H}_y-{\mathcal E}_{i_1}-{\mathcal E}_{i_2}-{\mathcal E}_{i_3}-{\mathcal E}_{i_4}\\
2{\mathcal H}_x+2{\mathcal H}_y-2{\mathcal E}_{i_1}-{\mathcal E}_{i_2}-{\mathcal E}_{i_3}-{\mathcal E}_{i_4}-{\mathcal E}_{i_5}\\
3{\mathcal H}_x+{\mathcal H}_y-{\mathcal E}_{i_1}-{\mathcal E}_{i_2}-{\mathcal E}_{i_3}-{\mathcal E}_{i_4}-{\mathcal E}_{i_5}-{\mathcal E}_{i_6}\\
3{\mathcal H}_x+2{\mathcal H}_y-2{\mathcal E}_{i_1}-2{\mathcal E}_{i_2}-{\mathcal E}_{i_3}-{\mathcal E}_{i_4}-{\mathcal E}_{i_5}-{\mathcal E}_{i_6}\\
4{\mathcal H}_x+{\mathcal H}_y-{\mathcal E}_{i_1}-{\mathcal E}_{i_2}-{\mathcal E}_{i_3}-{\mathcal E}_{i_4}-{\mathcal E}_{i_5}-{\mathcal E}_{i_6}-{\mathcal E}_{i_7}-{\mathcal E}_{i_8}
\end{array}
\end{align}  
and those with exchange of ${\mathcal H}_u$ and ${\mathcal H}_v$. 
Not all, but many pairs of these divisor classes admit a blow-down structure
for generic blow-up points.
For example, 
for
${\mathcal H}_u={\mathcal H}_x$ and ${\mathcal H}_u={\mathcal H}_x+{\mathcal H}_y-{\mathcal E}_{i_1}-{\mathcal E}_{i_2}$, 
${\mathcal F}_i$'s are given by
$${\mathcal H}_x-{\mathcal E}_{i_1},\ {\mathcal H}_x-{\mathcal E}_{i_2}\, 
{\mathcal E}_j (j\neq i_1,i_2) $$ 
and for 
${\mathcal H}_u={\mathcal H}_x+{\mathcal H}_y-{\mathcal E}_{i_1}-{\mathcal E}_{i_2}$ and ${\mathcal H}_u={\mathcal H}_x+{\mathcal H}_y-{\mathcal E}_{i_1}-{\mathcal E}_{i_3}$, 
${\mathcal F}_i$'s are given by
$${\mathcal H}_x-{\mathcal E}_{i_1},\ {\mathcal H}_y-{\mathcal E}_{i_1},\
 {\mathcal H}_x+{\mathcal H}_x-{\mathcal E}_{i_1}-{\mathcal E}_{i_2}-{\mathcal E}_{i_3},\
{\mathcal E}_j (j\neq i_1,i_2,i_3).$$ 
\end{example}

\begin{remark}
There is another way to obtain relatively minimal surface for elliptic surface case, though it needs heavy computation.   
Let $S$ be a rational elliptic surface (not necessarily minimal) where
the mapping $\varphi$ is lifted to an automorphism. 
Compute a $\R$-divisor $\theta$ by
$$\theta:=\lim_{n\to \infty} \frac{\tilde{\varphi}_*^n({\mathcal E})}{||\tilde{\varphi}_*^n({\mathcal E})||},$$
where $||\cdot||$ denotes the Euclidian norm of a divisor w.r.t. a fixed basis,    
and let $k>0$ be a minimum number such that $k\theta\in \pic(S)$.
Then, the linear system $|mk \theta|$ gives an elliptic fibration for
some integer $m\geq 1$ ($m$ is not always one,
(cf. Step 1 of Appendix of \cite{DF01} and the authors' paper \cite{CT12}).
Let $C$ be a curve in the linear system $|k\theta|$ 
(such $C$ exists \cite{CD89}).
By applying van Hoeji's algorithm \cite{Hoeij95} (cf. \cite{vgr}), we obtain a
birational transformation $S\to S', (x,y) \mapsto (u,v)$ such that $C$ is transformed into Weierstrass normal form $v^2=u^3-g_2u-g_3$.
Since the degree of this curve is three, 
$S'$ is obtained 
by $9$ blow-ups from $\P^2$. This implies $S'$ is a minimal elliptic
surface (the fibration is given by the linear system $|-mK_{S'}|$).
\end{remark}

\begin{remark}
If $\varphi$ is an automorphism of a non-minimal rational elliptic surface,
the invariant does not corresponds to the anti-canonical divisor,
because the self-intersection of the anti-canonical divisor
is negative in this case, while $\theta^2$ of the above remark should be zero.
\end{remark}

\subsection{A simple example which needs blowing down}
Let us show first a simple example which needs change of blow-down structure
to obtain relatively minimal surface. This example is due to Diller and Favre's paper \cite{DF01}:  (for simplicity we note $x_n=x, \bar{x}=x_{n+1}, \underline{x}=x_{n-1}$ and so forth
\begin{align}
&\left\{
\begin{array}{rcl}
\bar{x}&=&\displaystyle y+\frac{1}{2}\\[3mm]
\bar{y}&=&\displaystyle \frac{x(2y-1)}{2y+2}
\end{array}\right..
\end{align}
This system can be lifted to an automorphism on a surface $S$ by blowing up
 $\P^1\times \P^1$ at the singularity points of the dynamical systems:
\begin{align*}
& \E_1:(x,y)=(1,0),\
 \E_2 (1/2,-1/2),\
 \E_3 (0,-1),\
 \E_4 (-1/2,\infty),\ \\
&\E_5 (\infty,-1/2),\
 \E_6 (0,\infty),\
 \E_7 (\infty,0),\
 \E_8 (1/2,\infty),\
 \E_9 (\infty,1/2).
\end{align*}
Immediately one can see the action on the Picard group from the following singularity patterns:
\begin{align*}
&\H_y-\E_3\to\E_4\to\E_5\to \E_6\to \E_7\to\E_8\to\E_9\to \H_x-\E_1\\
&\H_y-\E_9\to \E_1\to \E_2\to \E_3\to \H_x-\E_4
\end{align*}
and also the invariant divisor classes
$\H_x+\H_y-\E_1-\E_2-\E_3$ and\\
$\H_x+\H_y-\E_4-\E_5-\E_6-\E_7-\E_8-\E_9$. The presence of invariant divisor calsses imposes making blow-down along the curve which corresponds to the divisor class $\H_x+\H_y-\E_1-\E_2-\E_3$ (it is the only one which has self-intersection -1, the other has self-intersection -3).
Hence we take the basis of blow-down structure as
\begin{align*}
&\H_u=\H_x+\H_y-\E_2-\E_3, \H_v=\H_x+\H_y-\E_1-\E_2,\\ 
&\H_x+\H_y-\E_1-\E_2-\E_3, 
\F_1=\H_x-\E_2,\ \F_2=\H_y-\E_2,\\
& \F_i=\E_{i+1}\ (i=3,4,5,6,7,8),
\end{align*}
where the linear systems of $\H_u$ and $\H_v$ are given by 
\begin{align*}
&|\H_u|: u_0(x-y-1)+u_1(2xy+x)=0,\\ 
&|\H_v|: v_0(x-y-1)+v_1(2xy-y)=0.
\end{align*}
Using these, we take the following change of variables:
\begin{align*}
u=\frac{2xy+x}{x-y-1},\quad v=\frac{2xy-y}{x-y-1},
\end{align*}
then our dynamical system (3) and (4) becomes
\begin{align}\left\{\begin{array}{rcl}
\bar{u}&=&\displaystyle \frac{2uv-u-v-1}{u-3v+1}\\
\bar{v}&=&\displaystyle \frac{-2uv}{u+v+1}
\end{array}\right..
\end{align}
This system has the following blow-up points:
\begin{align*}
&\F_1:(u,v)=(-1,0),\
\F_2(0,-1),\
\F_3(1,2),\
\F_4:(u,(v+1)/u)=(0,1).\\
&\F_5 (0,1),\
\F_6 (1,0),\
\F_7:((u+1)/v,v)=(1,0),\
\F_8 (2,1).
\end{align*}
and the linear system of the anti-canonical divisor gives
the invariant
$$K=\frac{u v(2u v - u - v - 1)}{(u - v)^2 - 1}=
\frac{x (2 x-1) y (2 y-1) (2x y-x+y+1)}{(x - y-1)^2}$$
and the invariant two form
$$\omega=\frac{du\wedge dv}{(u - v)^2 - 1}=\frac{dx\wedge dy}{1 - x + y}
.$$

\section{Discrete Nahm equations}
All the examples in this section preserve an elliptic fibration.

\subsection{Discrete Nahm equations with tetrahedral symmetry}
In \cite{PPS11}, Petrera, Pfadler and Suris proposed the following discretization of the reduced Nahm equations with tetrahedral symmetry
\begin{align}\label{Nahm1}&\left\{\begin{array}{rcl}
\bar{x}-x&=&\epsilon(x\bar x-y\bar y)\\
\bar{y}-y&=&-\epsilon(x\bar y+y\bar x)
\end{array}\right..
\end{align}
Here $\epsilon$ is related to the step of discretization. The integrability can be proved by the existence of the following conserved quantity and invariant two-form
\begin{equation}\label{inv3.1}
K=\frac{y(3x^2-y^2)}{-1+\epsilon^2(x^2+y^2)},\quad \omega=\frac{dx\wedge dy}{y(3x^2-y^2)}.
\end{equation}
In this case one can easily transform the system into a QRT one by the following variable transformation (about QRT mappings see \cite{qrt, Tsuda04, Duistermaat})
\begin{equation}\label{trans1}
u=\frac{1-\epsilon x}{y},v=\frac{1+\epsilon x}{y}.
\end{equation}
Immediately we get $\bar u=v$. From the equation \eqref{Nahm1} we get a QRT mapping 
$$3\bar{u}\underline{u}-u(\bar u+\underline{u})-u^2+4\epsilon^2=0$$
with the invariant:
$$K=\frac{-3(u-v)^2+4\epsilon^2}{2\epsilon^2(u+v)(u v-\epsilon^2)},\quad
\omega=\frac{du\wedge dv}{3(u - v)^2 - 4 \epsilon^2},$$
which are precisely \eqref{inv3.1} in the variables $x$ and $y$.

Now we are going to study the singularity structure and its space of initial conditions and recover the invariants. The fact that the conserved quantity is expressed by a ratio of a cubic polynomial implies that 
we have better to start with $\P^2$ than $\P^1\times \P^1$. 

On $\P^2:(X:Y:Z)=(x:y:1)$, we blow up the following points

\begin{align*}
&\E_1(-1:-\sqrt{3}: 2\epsilon),\
\E_2(1: \sqrt{3}:2\epsilon),\
\E_3(-1:\sqrt{3}:2\epsilon),\\
&\E_4(1:-\sqrt{3}:2\epsilon),\
\E_5(1:0:\epsilon),\
\E_6(-1:0:\epsilon),\\
&\E_7(1:0:0),\
\E_8(1:1:0),\
\E_9(1:-1:0).
\end{align*}

In order to blow down to $\P^1\times \P^1$, we take the basis of blow-down structure
$\H_x$, $\H_y$, $\F_1,\dots,\F_8$ as
\begin{align*}
&\H_x=\E-\E_5,\quad \H_y=\E-\E_6,\quad
\F_i=\E_i (i=1,2,3,4),\\ 
&\F_5=\E_7,\quad \F_6=\E_8,\quad \F_7=\E_8,\quad 
\F_8=E-\E_5-\E_7.
\end{align*}

The curves corresponding to the divisor classes $\H_x$ and $\H_y$ are:
$$\alpha_0(\epsilon X-Z)+\alpha_1 Y=0,\quad \beta_0(\epsilon X+Z)+\beta_1 Y=0.$$
They give immediately the change of variable
$$u=\frac{\epsilon x-1}{y},\quad v=\frac{\epsilon x+1}{y},$$
which is essentially \eqref{inv3.1} up to rescaling factors.

\subsection{Discrete Nahm equations with octahedral symmetry:}

The second Nahm equation is the one corresponding to octahedral symmetry. The system has the following form 
\begin{align}&\left\{\begin{array}{rcl}
\bar{x}-x&=&\epsilon(2x\bar x-12y\bar y)\\
\bar y-y&=&-\epsilon(3x\bar y+3y\bar x+4y\bar y)
\end{array}\right. ,
\end{align}
which is again integrable by the invariants:
\begin{align}\label{inv2-2}
K&=\frac{y(2x+3y)(x-y)^2}{1-10\epsilon^2(x^2+4y^2)+\epsilon^4(9x^4+272x^3y-352xy^3+696y^4)}\nonumber \\
\omega&=\frac{dx\wedge dy}{y (x - y) (2 x + 3 y)}.
\end{align}
Inspired by the transformation \eqref{trans1} we can simplify the system by the following transformations:
$$x=\frac{1}{3}(\chi-2y), \quad \bar x=\frac{1}{3}(\bar \chi-2\bar y)$$
and $u=(1-\epsilon\chi)/y, v=(1+\epsilon\chi)/y$. Finally we get a simpler equation but non-QRT type:
$$8\bar{u}\underline{u}-2u(\bar u+\underline{u})+20\epsilon(\bar u-\underline u)-4u^2+400\epsilon^2=0,$$
which can be written as a system on $\P^1\times\P^1$
\begin{align}\left\{\begin{array}{rcl}
\bar{u}&=v&\\
\bar{v}&=&\displaystyle \frac{(u+2v-20\epsilon)(v+10\epsilon)}{4u-v+10\epsilon}
\end{array}\right. .
\end{align}

The space of initial conditions is given by the $\P^1\times\P^1$ blown up at the following nine points:
\begin{align*}
&\E_1:(u,v)=(-10\epsilon,0),\
\E_2(0,10\epsilon),\
\E_3(10\epsilon,5\epsilon),\\
&\E_4(5\epsilon,0),\
\E_5(0,-5\epsilon),\
\E_6(-5\epsilon,-10\epsilon)\\
&\E_7(\infty,\infty),\
\E_8:(1/u,u/v)=(0,-1/2),\
\E_9:(1/u,u/v)=(0,-2).
\end{align*}

The action on the Picard group is the following:
\begin{align*}
&\bar{\H_u}=2\H_u+\H_v-\E_1-\E_3-\E_7-\E_8,\
\bar{\H_v}=\H_u\\
&\bar{\E_1}=\E_2,\
\bar{\E_2}=\H_u-\E_3,\
\bar{\E_3}=\E_4,\ \bar{\E_4}=\E_5,\ \bar{\E_5}=\E_6,\\
&\bar{\E_6}=\H_u-\E_1,\ 
\bar{\E_7}=\H_u-\E_8,\
\bar{\E_8}=\E_9,\
\bar{\E_9}=\H_u-\E_7.
\end{align*}
From this action one can see immediately that we have three invariant divisor classes:
\begin{align*}
&\alpha_0=\H_u+\H_v-\E_1-\E_2-\E_7,\ 
\alpha_1=\H_u+\H_v-\E_1-\E_2-\E_8-\E_9,\\
&\alpha_2=\E_7-\E_8-\E_9,\
\alpha_3=\H_u+\H_v-\E_3-\E_4-\E_5-\E_6-\E_7.
\end{align*}
The curve corresponding to $\alpha_0$ is a (-1) curve which must be blown down. Let $\H_a=\H_u+\H_v-\E_2-\E_7$ and $\H_b=\H_u+\H_v-\E_1-\E_7$, then their linear systems 
are given by
$$a_1 u+a_2(v-10\epsilon)=0,\quad b_1 (u+10\epsilon)+b_2 v=0$$
and the basis of blow-down structure is given by
\begin{align*}
&\H_a,\ \H_b,\ \alpha_0,\ \F_1=\H_u-\E_7,\ \F_2=\H_v-\E_7,\\ 
&\F_3=\E_3,\ \F_4=\E_4,\ \F_5=\E_5,\ \F_6=\E_6,\ \F_7=\E_8,\ \F_8=\E_9.
\end{align*}
So if we set:
$$a=\frac{v-10\epsilon}{u}\,\quad b=\frac{u+10\epsilon}{v},$$
our dynamical system becomes
\begin{align}\left\{ \begin{array}{rcl}
\bar{a}&=&\displaystyle \frac{3ab-2a+2}{a-4}\\
\bar{b}&=&\displaystyle \frac{4-a}{2a+1}
\end{array} \right..
\end{align}

This system has the following space of initial conditions which define a minimal rational elliptic surface:
\begin{align*}
&\F_1:(a,b)=(0,\infty),\
\F_2:(a,b)=(\infty,0),\\
&\F_3:(a,b)=(-1/2,4),\
\F_4:(a,b)=(-2,\infty)\\
&\F_5:(a,b)=(\infty,-2),\
\F_6:(a,b)=(4,-1/2),\\
&\F_7:(a,b)=(-2,-1/2),\
\F_8:(a,b)=(-1/2,-2).
\end{align*}
The invariants can be computed from the anti-canonical divisor as
$$K=\frac{(a b-1) (a b + 2a+2 b -5)}{4ab+2a+2b+1},\
\omega=\frac{da\wedge db}{(a b-1) (a b + 2a+2 b -5) } $$
which are equivalent to the invariants \eqref{inv2-2}.

\subsection{Discrete Nahm equations with icosahedral symmetry}

The last example of discrete reduced Nahm equations refers to icosahedral symmetry. It is given by 
\begin{align}\left\{\begin{array}{rcl}
\bar{x}-x&=&\epsilon(2x\bar x-y\bar y)\\
\bar y-y&=&-\epsilon(5x\bar y+5y\bar x-y\bar y)
\end{array}\right.
\end{align}
and is integrable as well. However the invariants here are more complicated. They are reported also by \cite{PPS11} as\footnote{
a sign in $c_2$ was corrected by information from the authors of that paper}
\begin{align}\label{Nahm3inv}
K=\frac{y(3x-y)^2(4x+y)^3}{1+\epsilon^2 c_2+\epsilon^4 c_4+\epsilon^6 c_6},
\quad 
\omega=\frac{dx\wedge dy}{y(3x-y)(4x+y)}
\end{align}
where
\begin{align*}
c_2&=-7(5x^2+y^2)\\
c_4&=7(37x^4+22x^2y^2-2xy^3+2y^4)\\
c_6&=-225x^6+3840x^5y+80xy^5-514x^3y^3-19x^4y^2-206x^2y^4.
\end{align*}

Again we can make first the following change of variable
$$x=\frac{1}{5}(X+\frac{y}{2}), \quad \bar x=\frac{1}{5}(\bar X+\frac{\bar y}{2}),$$
then we divide by $y\bar y$ both equations and call again $a=X/y, b=1/y, u=b-\epsilon a, v=b+\epsilon a$ and finally we get a simpler equation but non-QRT type:
$$6\bar{u}\underline{u}-u(\bar u+\underline{u})-\frac{7\epsilon}{2}(\bar u-\underline u)-4u^2+49\epsilon^2=0.$$

We can apply our procedure to this last non-QRT mapping. However, here we demonstrate that our procedure works well even 
for the original mapping.

The space of initial condition is given by the $\P^1\times\P^1$ blown up at the following 12 points:
\begin{align*}
&\E_1:(x,y)=(\infty, \infty),\
\E_2 (-1/7\epsilon, -3/7\epsilon),\
\E_3 (-1/7\epsilon, 4/7\epsilon),\\
&\E_4 (1/7\epsilon, 3/7\epsilon),\
\E_5 (1/7\epsilon, -4/7\epsilon)\,
\E_6 (1/5\epsilon,0),\\
&\E_7 (1/3\epsilon,0),\
\E_8 (1/\epsilon,0),\
\E_9 (-1/\epsilon,0),\\
&\E_{10} (-1/3\epsilon,0),\
\E_{11} (-1/5\epsilon,0),
\E_{12}:(1/x,x/y)=(0,1/3).
\end{align*}

On this surface the dynamical system is neither 
an automorphism nor analytically stable 
due to the following topological singularity patterns:
\begin{align*}
\H_y-\E_1\ (y=\infty)\to \mbox{point} \to \cdots \mbox{(4 points)} \cdots \to \mbox{point} \to \H_y-\E_1\\
\cdots \to \mbox{point} \to \mbox{point} \to \H_x-\E_1\ (x=\infty) \to
\mbox{point} \to \mbox{point} \to \cdots,
\end{align*}
where the image of a curve under $\varphi^n$ is defined as 
(ia) in Section 1.
Moreover, the curve $4x+y=0: \H_x+\H_y-\E_1-\E_3-\E_5$ is invariant.
We blow down along these three curves with the blow-down
structure
\begin{align*}
&\H_u=\H_x+\H_y-\E_1-\E_3,\ \H_v=\H_x+\H_y-\E_1-\E_5,\\ 
&\H_x-\E_1,\ \H_y-\E_1,\ \H_x+\H_y-\E_1-\E_3-\E_5,\\
&\F_1=\E_{12},\ \F_2=\E_2,\ \F_3=\E_4,\ \F_4=\E_6,\\ 
&\F_5=\E_7,\ \F_6=\E_8,\
\F_7=\E_9,\ \F_8=\E_{10},\ \F_9=\E_{11},
\end{align*}
where the linear systems of $\H_v$ and $\H_v$ are given by
\begin{align*}
|\H_u|:& u_0(1+7 \epsilon x)+u_1(4x+y)\\
|\H_v|:& v_0(1-7 \epsilon x)+v_1(4x+y).
\end{align*}
If we take the new variables $u$ and $v$ as
$$v=\frac{2(1+7 \epsilon x)}{\epsilon(4x+y)},\
v=\frac{2(1-7 \epsilon x)}{\epsilon(4x+y)},
$$
then we have
\begin{align*}
&\F_1:(u,v)=(2,-2), \F_2:(0,-4),
\F_3:(4,0),
\F_4:(6, -1), \F_5:(5, -2),\\
&\F_6:(4, -3), \F_7:(3, -4), \F_8:(2, -5), \F_9:(1, -6).
\end{align*}
The dynamical system becomes an automorphism having
the following topological singularity patterns
\begin{align*}
&\H_v-\F_9\to \F_2\to \F_1 \to \F_3\to \H_u-\F_4\\
&\H_v-\F_3\to \F_4 \to \F_5 \to \F_6 \to \F_7 \to \F_8 \to \F_9 \to \H_u-\F_2 
\end{align*}
and $\H_u\to \H_u+\H_v-\F_2-\F_4$. Hence we find
the invariant $(-1)$ curve $\H_u+\H_v-\F_1-\F_2-\F_3$, which should be blown down.
Again we take the blow-down structure as
\begin{align*}
&\H_s=\H_u+\H_v-\F_1-\F_2,\ \H_t=\H_u+\H_v-\F_1-\F_3,\\ 
&\H_u+\H_v-\F_1-\F_2-\F_3,\ \F_1'=\H_a-\F_1,\ \F_2'=\H_b-\F_1 \\
&\F_3'=\F_4,\ \F_4'=\F_5,\ \F_5'=\F_6,\ \F_6'=\F_7,\\ 
&\F_7'=\F_8,\ \F_8'=\F_9,
\end{align*}
where the linear systems of $\H_s$ and $\H_t$ are given by
\begin{align*}
|\H_s|:& s_0u(v+2)+s_1(u-v-4)\\
|\H_t|:& t_0v(u-2)+t_1(u-v-4) 
\end{align*}
and hence we take the new variables $s$ and $t$ as
$$s=-\frac{3u(v+2)}{2(u-v-4)},\
t=-\frac{3v(u-2)}{2(u-v-4)}.
$$
Then we have
\begin{align*}
&\F_1':(s,t)=(3,0),\ \F_2'(0,3),\
\F_3'(-3,2),\
\F_4':(\frac{s}{t-3},d-3)=(5,0),\\ 
&\F_5'(2, 3),\ \F_6'(3, 2),\ \F_7':(u-3,\frac{t}{s-3})=(0,5),\ \F_8'(2, -3)
\end{align*}
and
\begin{align*}&\left\{ \begin{array}{rcl}
\bar s&=&\displaystyle \frac{2s t -3s-3t+9}{s+t-3}\\[2mm] 
\bar t&=&\displaystyle \frac{2(s-3)(t+3)}{3s-t-9}
\end{array}\right..
\end{align*}
The invariants can be computed by using the 
the anticanonical divisor as
\begin{align} K'=\frac{(s-t)^2+4(s+t)-21}{(s-2)(t-2)(2st-5s-5t+15)} 
=\frac{-56 \epsilon^6 y (-3x + y)^2 (4x + y)^3}{d_1d_2d_3}
\end{align}
and 
\begin{align}
\omega=\frac{2\epsilon ds\wedge dt}{(s - t)^2 + 4(s + t) - 21}
=\frac{dx\wedge dy}{y(3x-y)(4x+y)},
\end{align}
where
\begin{align*} 
d_1&= -3 - 12\epsilon x + 15 \epsilon^2 x^2 - 3 \epsilon y 
- 17 \epsilon^2 x y + 4\epsilon^2y^2\\
d_2&= -3 + 12 \epsilon x + 15 \epsilon^2 x^2 + 3 \epsilon y 
- 17 \epsilon^2 x y + 4 \epsilon^2 y^2\\
d_3&= -3 + 27 \epsilon^2 x^2 + 10\epsilon^2 x y + 10 \epsilon^2 y^2.
\end{align*}
The denominator of $K'$ is related to $K$ of \eqref{Nahm3inv} as 
$$d_1d_2d_3=160\epsilon^6  (\mbox{numerator of }K)
-27 (\mbox{denominator of }K).$$

\section{Linerizable mappings}
In this section we demonstrate that our method works well 
also for linearizable mappings.
The first example is a simple 
 non-autonomous linearizable mapping studied in \cite{TEGORS}.
We show our method is different from that paper and \cite{DF01}. 
The second example is also a linearizable mapping 
proposed again by \cite{PPS11} as a discretization of the Suslov system.

\subsection{A non-autonomous linearizable mapping}

Here we consider the following very simple mapping 
\begin{align}&\left\{\begin{array}{rcl}
\bar{x}&=&y\\
\bar {y}&=&\left(-\frac{y}{x}+a_n \right)y
\end{array}\right.,
\end{align}
where $a_n$ is an arbitrary sequence of complex numbers.
This dynamical system is a linearizable mapping studied in \cite{TEGORS}
and the degree of this dynamical system grows linearly and
it is lifted to an analytically stable mapping by blowing up
at the following points:
\begin{align*}
&\E_1:(x,y)=(0,0),\ 
\E_2:(\infty,\infty).
\end{align*}
The topological singularity patterns are
\begin{align*}
&(\frac{x}{y},y)=(0,0)\to \H_x-\E_1 \to \H_y-\E_2 \to (\frac{1}{x},\frac{x}{y})=(0,0) \\&\mbox{(point on $\E_2$)} \to \H_x-\E_2 \to \mbox{(curve)}\\
&\mbox{(curve)} \to \H_y-\E_1 \to \mbox{(point on $\E_2$)}.
\end{align*}
These are not confined at all.
Moreover, we can compute the action on the Picard group as
\begin{align*}
&\bar \H_x=2\H_x+\H_y-\E_1-\E_2\\
&\bar \H_y=\H_x,\quad \bar \E_1=\H_x,\quad \bar \E_2=\H_x.
\end{align*}
However, since the dynamical system is not an isomorphism, 
we need to compute very carefully for this result.  
One can see detail of such computation in \cite{TEGORS}. 

In this paper, we shall linearize the dynamical system with 
using singularity patterns
instead of the action on the Picard group.

From the singularity pattern, we can blow down the surface 
along $\H_x-\E_1$, keeping analytical stability. 
Then we can easily find a basis of blow-down structure as
$$\H_u=\H_x,\ \H_v=\H_x+\H_y-\E_1-\E_2,\ \F_1=\H_x-\E_1,\ \F_2=\H_x-\E_2.$$
where the linear systems of $\H_u$ and $\H_v$ are
$$|\H_u|: u_0x+u_1=0,\quad |\H_v|: v_0x+v_1 y=0.$$
Taking new variables $u$ and $v$ as
$u=x$ and $v=y/x$, we have
\begin{align}&\left\{\begin{array}{rcl}
\bar{u}&=&uv\\
\bar {v}&=& v+a_n
\end{array}\right..
\end{align}

\subsection{Discrete Suslov system}

The discrete Suslov system proposed in \cite{PPS11} is a linearizable mapping:
\begin{align}&\left\{\begin{array}{rcl}
\bar{x}-x&=&\epsilon a(\bar x y+ x \bar y)\\
\bar y-y&=&-2\epsilon x\bar x
\end{array}\right. .
\end{align}
Again, the degree of this dynamical system grows linearly and
it is lifted to an analytically stable mapping by blowing up
at the following points: (we put $a=-b^2$ for simplicity)
\begin{align*}
&\E_1:(x,y)=\left(-\frac{1}{b\epsilon},\frac{1}{b^2 \epsilon}\right),\
\E_2:\left(\frac{1}{b\epsilon},\frac{1}{b^2 \epsilon}\right),\\
&\E_3:\left(-\frac{1}{b\epsilon},-\frac{1}{b^2 \epsilon}\right),\
\E_4:\left(\frac{1}{b\epsilon},-\frac{1}{b^2 \epsilon}\right),\ 
\E_5: (\infty,\infty).
\end{align*}
The topological singularity patterns are
\begin{align*}
&x=\infty \to (0,-\frac{1}{b^2\epsilon})\\
&y=\infty \to y=\infty\\
&(2bex+b^2\epsilon y+1=0)\to \E_3 \to \E_2\to (2b\epsilon x-b^2\epsilon y+1=0)\\
&(-2b\epsilon x+b^2\epsilon y+1=0)\to \E_4 \to \E_1\to (-2b\epsilon x-b^2\epsilon  y+1=0)\\
&(2b^2 \epsilon^2 x^2- b^2 \epsilon y-1=0)\to \E_5\to (2b^2 \epsilon^2 x^2+ b^2 \epsilon y-1=0),
\end{align*}
where divisor classes are
\begin{align*}
&x=\infty : \H_x-\E_5\\
&x=\infty : \H_y-\E_5\\ 
&2b\epsilon x-b^2\epsilon y+1=0 : \H_x+\H_y-\E_4-\E_5 \\
&-2b\epsilon x-b^2\epsilon  y+1=0 : \H_x+\H_y-\E_3-\E_5 \\
&2b^2 \epsilon^2 x^2+ b^2 \epsilon y-1=0 : 2\H_x+\H_y-\E_3-\E_4-\E_5.
\end{align*}

At first, we blow down along $\H_x-\E_5$ and $\H_y-\E_5$.
For that purpose we take the blow-down structure as
\begin{align*}
&\H_s:=\H_x+\H_y-\E_1-\E_5,\ \H_t:=\H_x+\H_y-\E_2-\E_5,\\
&\H_x-\E_5,\ \H_y-\E_5,\ \H_x+\H_y-\E_1-\E_2-\E_5,\ \E_3,\ \E_4.
\end{align*}
Then we have a surface whose Picard group is generated by
$\H_s$, $\H_t$, $\E_3$, $\E_4$ where the dynamical system
is still analytically stable. We abbreviate detail, but
again we find effective (-1) divisor classes $\H_s-\E_3$ and $\H_s-\E_4$ 
in singularity pattern which can be blown down preserving
analytical stability. 
Hence we take a basis of blow-down structure as
\begin{align*}
&\H_u:=\H_s+\H_t-\E_3-\E_4=2\H_x+2\H_y-\E_1-\E_2-\E_3-\E_4-2\E_5,\\
&\qquad u_0(x^2-b^2 y^2)+u_1(1-b^2\epsilon^2 x^2)=0,\\
&\H_v:=\H_s=\H_x+\H_y-\E_1-\E_5: v_0(1+b\epsilon x)+v_1(x+by)=0\\
&\H_s-\E_3=\H_x+\H_y-\E_1-\E_3-\E_5\\
&\H_s-\E_4=\H_x+\H_y-\E_1-\E_4-\E_5.
\end{align*}
If we take the new variables $u$ and $v$ as
$$ u = \frac{x^2 - b^2 y^2}{1 - b^2 \epsilon^2 x^2},\quad
 v = \frac{1 + b \epsilon x}{x + b y},$$
then the dynamical system becomes
\begin{align}&\left\{\begin{array}{rcl}
\bar u&=&u\\
\bar v&=&\displaystyle \frac{b \epsilon+v}{1-b \epsilon u v}
\end{array}\right. .
\end{align}

\begin{remark}
The action of the mapping on the Picard group on the first
surface is given by 
\begin{align*}
&\bar \H_x = 2\H_x+\H_y- \E_3-\E_4-\E_5\\
&\bar \H_y = 2\H_x+2\H_y- \E_3-\E_4-2\E_5\\
&\bar \E_1 = 2\H_x+\H_y- \E_3-\E_5\\
&\bar \E_2 = 2\H_x+\H_y- \E_4-\E_5\\
&\bar \E_3 = \E_2,\quad \bar \E_4= \E_1\\
&\bar \E_5 = 2\H_x+\H_y- \E_3-\E_4-\E_5
\end{align*}
and $\H_u$ is the invariant divisor class whose self-intersection is zero. 
\end{remark}

\noindent{\bf Acknowledgements:} A.S.Carstea was supported by the project IDEI, PN-II-ID-PCE-2011-3-0137, Romanian Ministery of Education.

\end{document}